\documentclass{amsart}

\usepackage{amssymb}
\usepackage{amsmath}
\usepackage{amsthm}
\usepackage[mathscr]{eucal}
\usepackage{geometry,xcolor,graphicx}
\usepackage{mathpazo}

\newtheorem*{main*}{Main Theorem}

\theoremstyle{definition}

\theoremstyle{remark}

\numberwithin{equation}{section}

\newcommand{\dbar}{$\bar{\partial}$}
\newcommand{\mdbar}{\bar{\partial}}

\usepackage{mathrsfs,latexsym,url,setspace}
\onehalfspace

\begin{document}

\bibliographystyle{plain}

\title[Regularity of projections and of
 \dbar]{The regularity of projection operators and solution operators to
\dbar\
 on weakly pseudoconvex domains}

\author[D. Ehsani]
 {Dariush Ehsani}
\thanks{This research was (partially) supported by the Deutsche Forschungsgemeinschaft (DFG, German Research Foundation),
grant RU 1474/2 within DFG�s Emmy Noether Programme.}

\address{
Hochschule Merseburg\\
Fachbereich IKS\\
Geusaer Stra\ss e\\
 D-06217 Merseburg\\
 Germany}
 \email{dehsani.math@gmail.com}


\begin{abstract}
 We relate the existence and regularity of a solution operator to
 $\mdbar$ on smoothly bounded pseudoconvex domains
 to the existence and regularity of a projection operator
 onto the kernel of $\mdbar$.
\end{abstract}

\maketitle

\section{Introduction}

  In this article we show that the
 existence of a solution operator to the \dbar-equation
  satisfying estimates on Sobolev spaces of smoothly bounded pseudoconvex domains
  can be deduced from the existence of a continuous projection operator.
  We let $\Omega\subset\mathbb{C}^n$ be a smoothly bounded pseudoconvex domain.
   We denote the Sobolev spaces by
 $W^s_{(p,q)}(\Omega)$, the space of
 $(p,q)$-forms whose derivatives of order $\le s$ (of its
 coefficients)
are in  $L^2_{(p,q)}(\Omega)$.   Let $D_1,\ldots {,}D_{2n}$ be an
orthonormal frame of $2n$ vector fields.  The norm $\|\cdot\|_s$
attached to $W^s_{(p,q)}(\Omega)$ is given by
\begin{equation*}
 \|f\|_s = \sum_{k=0}^s \sum_{j=1}^{2n}
 \|D_j^kf\|_{L^2_{(p,q)}(\Omega)}{.}
\end{equation*}
 We look at the
 relation between the existence of a projection operator,
  ${\bf T}_{(p,q)} :L^2_{(p,q)}(\Omega)\rightarrow L^2_{(p,q)}(\Omega)\cap \mbox{ker}(\mdbar)$,
   which is also
  a linear operator between Sobolev spaces, so that
${\bf T}_{(p,q)} :W^s_{(p,q)}(\Omega)\rightarrow
W^s_{(p,q)}(\Omega)$ for all $s\ge 0$
   and the
  existence of a solution operator,
 ${\bf S}_{(p,q+1)} :L^2_{(p,q+1)}(\Omega)\cap \mbox{ker}(\mdbar) \rightarrow L^2_{(p,q)}(\Omega)$
  such that $\mdbar{\bf S}_{(p,q+1)} f = f$
  and $\|{\bf S}_{(p,q+1)} f\|_s\lesssim \|f\|_s$ for all {$s\ge 0$, where
  the constant in the inequality is independent of $f$}.  Our Main Theorem
 in this regard is
 \begin{main*}
  Let $\Omega\subset \mathbb{C}^n$ be a smoothly bounded psedoconvex domain.
   For $0\le p\le n$ and $0\le q\le n-1$,
 suppose there exists a projection operator
 ${\bf T}_{(p,q)} :L^2_{(p,q)}(\Omega)\rightarrow L^2_{(p,q)}(\Omega)\cap \mbox{ker}(\mdbar)$,
 {with the property that ${\bf T}_{(p,q)}$ is continuous as an
 operator}
\begin{equation*}
{\bf T}_{(p,q)} :W^s_{(p,q)}(\Omega)\rightarrow
W^s_{(p,q)}(\Omega) \qquad \forall s\ge 0.
\end{equation*}
  Then there exists a solution operator, ${\bf S}_{(p,q+1)}$, such
  that
$\mdbar \circ {\bf S}_{(p,q+1)} = I$ on \dbar-closed
$(p,q+1)$-forms
 which is continuous as {an operator}
\begin{equation*}
 {\bf S}_{(p,q+1)}: W^s_{(p,q+1)}(\Omega)\cap
\mbox{ker}(\mdbar) \rightarrow
 W^s_{(p,q)}(\Omega) {\qquad \forall s\ge 0}.
\end{equation*}
 \end{main*}

{Naively, one could also obtain estimates for a projection
operator given a solution operator: given ${\bf S}_{(p,q+1)}$,
then we could form the projection ${\bf T}_{(p,q)}= I - {\bf
S}_{(p,q+1)}\circ \mdbar$.  Due to the \dbar-operator, however,
there is a loss of one derivative in the reverse direction of the
theorem.  One would need estimates for the composite operator
${\bf S}_{(p,q+1)}\circ \mdbar$ to obtain a theorem with the
reverse implication.  Such an equivalence between the regularity
of the Bergman projection and the \dbar-Neumann operator (and
through it the regularity of the canonical solution to the
\dbar-problem) was proved in \cite{BSt90}.  Thus, if the Bergman
projection $P_{(p,q)}$ on $(p,q)$-forms satisfies the hypothesis
of the theorem, we can take ${\bf S}_{(p,q+1)}={\bf C}_{(p,q+1)}$,
the canonical solution operator, or the operator giving the
solution of minimal $L^2$-norm.
 However,}
 we note that the Bergman projection does not always hold the
 desired properties of the projection operators in the Main Theorem.
 It is a result of
 Barrett that there exist smoothly bounded pseudoconvex domains on
 which the Bergman projection fails to map $W^s(\Omega)$ to itself
 for large enough $s$ \cite{Ba92}.  Furthermore, the solution
 operators to the \dbar-equation on weakly pseudoconvex domains of Kohn
  \cite{K73} do not satisfy the desired properties of the solution operators
 of the Main Theorem.  The solution operator
  of Kohn is only shown to give regularity up to a fixed level;
   for a fixed $s$ a solution operator, {$K_{s,q+1}$,}
 can be constructed such that for a \dbar-closed $(0,q+1)$-form,
 $f$ we have $\mdbar {K_{s,q+1}} f = f$ and ${K_{s,q+1}}:W^k_{(0,q+1)}(\Omega)\rightarrow
 W^k(\Omega)$ for $k\le s$.  We are interested in the estimates holding simultaneously for all $s$.

The Main Theorem is proved by construction.  Given a projection, a
solution operator is constructed.  Thus, the Main Theorem suggests
one method of obtaining (regular) solution operators to \dbar.
Namely, if one can find a projection operator which preserves
Sobolev spaces, then a solution operator to the \dbar-equation can
be found with no loss of derivatives.

\section{Proof of the Main Theorem}

It suffices to work with
 $(0,q)$-forms, and we simply drop the $p$ component in the
 indices, e.g. $W^s_{q}(\Omega):=W^s_{(0,q)}(\Omega)$, etc.

Let ${\bf C}_{q+1}$ be the canonical solution operator (i.e. the
 operator which gives the solution of minimal $L^2$-norm)
 to the \dbar-equation, and let ${\bf P}_q$ denote the
  Bergman projection.
We note the relation
\begin{equation}
 \label{dbn}
{\bf P}_q={\bf I}- {\bf C}_{q+1}\circ\mdbar.
\end{equation}
We also denote by $K_{s,q+1}$ the solution of Kohn which maps
$W^s_{q+1}(\Omega)$ to $W^s_{q}(\Omega)$.

We show how to obtain a solution operator ${\bf S}_{q+1}$ by
combining the operators ${\bf C}_{q+1}$, ${\bf P}_q$, and ${\bf
T}_{q}$.  We thus assume the existence of ${\bf T}_{q}$ as in the
Main Theorem, and we define
\begin{equation}
 \label{level0}
{\bf S}_{q+1}= {\bf C}_{q+1}+({\bf P}_q-{\bf T}_q)\circ K_{0,q},
\end{equation}
restricted to $\mbox{ker}(\mdbar)$.  We note that for $s_1,s_2\ge
0$, $K_{s_1,q+1}-K_{s_2,q+1}$ maps \dbar-closed $(0,q+1)$-forms
onto \dbar-closed $(0,q)$-forms.  And, as both ${\bf P}_q$
 and ${\bf T}_q$ reproduce \dbar-closed $(0,q)$-forms, we have
\begin{equation*}
({\bf P}_q-{\bf T}_q) \circ K_{s_1,q+1}\equiv  ({\bf P}_q-{\bf
T}_q)\circ K_{s_2,q+1}
\end{equation*}
   on
 $W^{\max\{s_1,s_2\}}_{q+1}(\Omega) \cap
\mbox{ker} (\mdbar)$.  We thus conclude that
 with ${\bf S}_{q+1}$ defined as in \eqref{level0} for any $s\ge{0}$,
 we also have
\begin{equation*}
 {\bf S}_{q+1}= {\bf C}_{q+1}+({\bf P}_q-{\bf T}_q)\circ K_{s,q},
\end{equation*}
 on $W^{\infty}_{q}(\Omega) \cap \mbox{ker} (\mdbar)$ .

 Because $\mdbar \circ({\bf P}_q-{\bf T}_q) \equiv 0$,
 for $f\in L^2_{q+1}(\Omega)$ a \dbar-closed form we have
\begin{align*}
 \mdbar {\bf S}_{q+1} f &= \mdbar {\bf C}_{q+1} f +
 \mdbar ({\bf P}_q-{\bf T}_q) K_{0,q} f\\
 &=f.
 \end{align*}

We finish the proof of the Main Theorem by showing that for all
$s\ge0$ we have
  $ \| {\bf S}_{q+1}f \|_s  \lesssim \|f\|_s$  if $f\in W^s_{q+1}(\Omega)\cap \mbox{ker} (\mdbar)$.
    With $f\in W^s_{q+1}(\Omega)\cap \mbox{ker} (\mdbar)$ we have
\begin{align*}
 {\bf S}_{q+1}f
 &={\bf C}_{q+1}  f +
  ({\bf P}_q-{\bf T}_q)  K_{0,q} f\\
 &={\bf C}_{q+1} \mdbar K_{s,q}  f +
  ({\bf P}_q-{\bf T}_q)  K_{s,q} f\\
 &= K_{s,q} f-{\bf T}_q K_{s,q} f,
\end{align*}
 by \eqref{dbn}.
 By assumption we have that
 ${\bf T}_q$ preserves $W^s_q(\Omega)$, and we
conclude
\begin{align*}
\| {\bf S}_{q+1}f \|_s & \lesssim \| K_{s,q+1} f \|_s\\
 & \lesssim \|f\|_s.
\end{align*}
This finishes the proof of the Main Theorem.


\begin{thebibliography}{1}

\bibitem{Ba92}
D.~Barrett.
\newblock Behavior of the {Bergmann} projection on the {Diederich-Forn\ae ss}
  worm.
\newblock {\em Acta. Math.}, 168:1--10, 1992.

\bibitem{BSt90}
Boas, H. and Straube, E.
\newblock Equivalence of regularity for the {Bergman} projection and the
  $\bar{\partial}$-{Neumann} operator.
\newblock {\em Manuscripta Math.}, 67:25--33, 1990.


\bibitem{K73}
J.J. Kohn.
\newblock Global regularity for $\bar{\partial}$ on weakly pseudo-convex
  manifolds.
\newblock {\em Trans. Amer. Math. Soc.}, 181:273--292, 1973.

\end{thebibliography}
\end{document}